\documentclass[11pt]{amsart}
\usepackage{ifthen,amssymb,graphicx}
\usepackage[all]{xy}

\newcommand{\dbc}[1]{{D}^b(#1)}

\newcommand{\fmf}[3]{{\Phi^{#1}_{{\scriptscriptstyle #2\!\rightarrow\! #3}}}}

\newcommand{\bR}{{\mathbf R}}
\newcommand{\bL}{{\mathbf L}}

\newcommand{\bZ}{\mathbb{Z}}

\newcommand{\cE}{{\mathcal E}}
\newcommand{\cF}{{\mathcal F}}
\newcommand{\cG}{{\mathcal G}}

\newcommand{\cI}{{\mathcal I}}

\newcommand{\calL}{{\mathcal L}}
\newcommand{\cK}{{\mathcal K}}
\newcommand{\cM}{{\mathcal M}}
\newcommand{\cN}{{\mathcal N}}
\newcommand{\cO}{{\mathcal O}}
\newcommand{\cP}{{\mathcal P}}
\newcommand{\cQ}{{\mathcal Q}}

\newcommand{\cS}{{\mathcal S}}
\newcommand{\cT}{{\mathcal T}}
\newcommand{\cU}{{\mathcal U}}

\newcommand{\PP}{B}

\newcommand{\f}{\mathfrak{f}}

\newcommand{\lotimes}{{\,\stackrel{\mathbf L}{\otimes}\,}}

\newcommand{\rest}[2]{{#1}_{\vert #2}}
\newcommand{\what}[1]{{\widehat #1}}

\newcommand\rk{\operatorname{rk}}
\newcommand\ch{\operatorname{ch}}

\newcommand\Hom{\operatorname{Hom}}

\newcommand\Ext{\operatorname{Ext}}
\newcommand{\SHom}{{\mathcal{H}om}}

\newcommand{\cplx}[1]{\mathcal{#1}^{\scriptscriptstyle\bullet}}

\newcommand\iso{\kern.35em{\raise3pt\hbox
{$\sim$}\kern-1.1em\to}\kern.3em}

\newcommand{\marginnote}[1]{\ifthenelse{\isodd{\thepage}}{\normalmarginpar}
{\reversemarginpar}\marginpar{\fbox{\parbox{24mm}{\sloppy\footnotesize #1}}}}

\newtheorem{thm}{Theorem}[section]
\newtheorem{lem}[thm]{Lemma}
\newtheorem{cor}[thm]{Corollary}
\newtheorem{prop}[thm]{Proposition}

\theoremstyle{definition}
\newtheorem{rem}[thm]{Remark}
\newtheorem{defin}[thm]{Definition}
\newtheorem{exa}[thm]{Example}

\def\hra{\hookrightarrow}

\newcommand\Td{\operatorname{Td}}

\title{STABLE SHEAVES OVER K3 FIBRATIONS}

\author[B. Andreas]{Bj\"orn Andreas}
\address{  Departamento de Matem\'aticas, Universidad de Salamanca, Plaza
de la Merced 1-4, 37008 Salamanca, Spain}
\email{bandreas@usal.es}

\author[D. Hern\'andez Ruip\'erez]{Daniel Hern\'andez Ruip\'erez}
\author[D. S\'anchez G\'omez]{Dar\'{\i}o S\'anchez G\'omez}
\address{Departamento de Matem\'aticas {\rm and }Instituto Universitario de F\'{\i}sica Fundamental y Matem\'aticas
(IUFFYM), Universidad de Salamanca, Calle del Parque s/n, 37008 Salamanca, Spain}
\email{ruiperez@usal.es dario@usal.es}

\date{\today}
\begin{document}
\begin{abstract}
We construct stable sheaves over K3 fibrations using a relative Fourier-Mukai transform
which describes the sheaves in terms of spectral data similar to the construction for elliptic
fibrations. On K3 fibered Calabi-Yau threefolds we show that the Fourier-Mukai transform induces an embedding ion of the relative Jacobian of spectral line bundles on spectral covers into the moduli space of sheaves of given invariants. This makes the moduli space of spectral sheaves to a generic torus fibration over the moduli space of curves of given arithmetic genus on the Calabi-Yau manifold.
\end{abstract}

\maketitle
\section{Introduction}\label{intro}
Moduli spaces of stable sheaves and bundles have been extensively studied
on complex curves and surfaces but on higher dimensional varieties not much is known
about them. Aside from their mathematical importance, these moduli spaces provide a
geometric background for various aspects in string theory such as heterotic string compactifications
or the description of D-branes on complex three-dimensional Calabi-Yau varieties, which occur
in the mirror symmetry program.

If the variety admits a fibration structure, a natural procedure to describe moduli
spaces of sheaves or bundles is to first construct them fiberwise and then to find an appropriate
global description. This method has been successfully employed by Friedman, Morgan and Witten \cite{FMW99} to construct stable vector bundles on elliptic fibrations in terms of spectral covers. Equivalently, this construction can be understood as a relative Fourier-Mukai transform. More precisely, by the invertibility of the Fourier-Mukai transform one gets for elliptic fibrations
a one-to-one correspondence between fiberwise torsion free semistable sheaves of rank $n$ and degree 0 and torsion free rank one sheaves on spectral covers.

This paper has been motivated by the question whether K3 fibrations admit a similar description
of stable sheaves and their moduli spaces in terms of spectral data, and as such the paper continues  earlier studies of Thomas \cite{Th00a} as well as of Bridgeland and Maciocia \cite{BrM02}. A physical motivation of the paper is the search for stable bundles and sheaves which provide new
solutions to heterotic string theory, or describe D-brane configurations on Calabi-Yau spaces. In particular, the class of $K3$-fibrations has been studied extensively in the physical literature
whereas the existence of stable sheaves was assumed but no explicit constructions given (cf.~\cite{kllw}).

By analogy to the elliptic fibrations it is suggestive to first recall the description of sheaves
on a single K3 surface and the corresponding dual variety (the analog of the generalized Jacobian
parametrizing torsion free rank one sheaves on an elliptic curve). For this consider a smooth
algebraic K3 surface $S$ with a fixed ample line bundle $H$ on it. The role of the dual
variety is played by the moduli space $\cM^H(v)$ of Gieseker-stable sheaves $\cF$ on $S$
having a fixed Mukai vector
$$
v(\cF)=\ch(\cF)\sqrt{\Td(S)}=(r,l,s)\in H^0(S,{\mathbb Z})\oplus H^{1,1}(S,{\mathbb Z})\oplus H^4(S, {\mathbb Z})\,.
$$
The moduli space $\cM^H(v)$ is a quasi-projective subscheme of the moduli space of simple sheaves on $S$. A
well-known result of Mukai states that if $\cM^H(v)$ is nonempty and compact of dimension two, then
it is a K3 surface isogenous to the original K3 surface $S$. Moreover, if in addition the Mukai vector
is primitive and satisfies $\gcd(r, l H, s)=1$ then $\cM^H(v)$ is a fine moduli space
and parametrizes stable sheaves with Mukai vector $v$ so that there is a universal sheaf on
$S\times \cM^H(v)$.

One approach to globalize the construction is as follows (for a different approach cf.~\cite{Th00a}). Let $p\colon X\to B$ be a K3 fibration and $Y$ a component of the relative moduli space of stable
sheaves (with generic fiber $\cM^H(v)$) on the fibers of $p$, which is a fine moduli space, so that there is a universal sheaf $\cP$ on $Y\times X$.  If, in addition, $X$ and $Y$ have the same dimension and ${\hat p}\colon Y\to B$ is equidimensional, Bridgeland and Maciocia \cite{BrM02} proved that $Y$ is a non-singular projective variety, ${\hat p}\colon Y\to B$ is a K3 fibration and the integral functor $\dbc{Y}\to \dbc{X}$ with kernel $\cP$ is an equivalence of derived categories, that is, a Fourier-Mukai transform. A more detailed review of the relative moduli space construction and Fourier-Mukai transform is given in Section \ref{section:moduli}. 

In Section \ref{section:stable} we construct torsion free sheaves on a K3 fibration $X$ which are $\mu$-stable with respect to some appropriate polarization. For the construction of sheaves on $X$ we are interested in transforming sheaves into sheaves rather than complexes using the Fourier-Mukai transform of \cite{BrM02}. For this it is useful to adopt the notion of ${\rm WIT}_i$ (Weak Index Theorem) sheaves.
We say an $\cO_X$-module $\cE$ is ${\rm WIT}_i$ if its Fourier-Mukai transform is concentrated in degree $i$.
We introduce an appropriate class of sheaves which satisfy the ${\rm WIT}_i$ condition. In analogy to the elliptic fibration set-up, we consider a $n$-fold cover $C$ (the {\it spectral cover}) over $B$ and prove that whenever $C$ is integral and $\calL$ a line bundle over $C$ (which is shown to be $WIT_0$), then the Fourier-Mukai transform of ${\calL}$ leads to a single torsion free sheaf $\widehat\cE$ (the {\it spectral sheaf}) on $X$. We prove that $\widehat\cE$ is $\mu$-stable with respect to $H+M\f$ with $M>0$ and $H$ a fixed ample divisor in $X$ and $\f$ the class of the fiber of $p$ (Theorem \ref{thm:mu-stability}).
As a result we get $\mu$-stable torsion free sheaves under weaker assumptions than in \cite{Th00a}.
Under the assumptions of \cite{Th00a} we prove that the cover $C$ can be chosen such that
the spectral sheaf $\widehat\cE$ is locally free.
In Subsection \ref{topinv} we restrict $X$ to be a K3 fibered Calabi-Yau threefold, that is, throughout this section we will assume that $X$ has a trivial canonical bundle, and compute under this assumption the topological invariants of $\widehat\cE$.

In Section \ref{section:ext} we prove that a non-split extension of stable spectral bundles on a K3 fibration is stable with respect to $H+M\f$, provided the slopes of the input sheaves satisfy certain conditions, which are required to rule out possible destabilizing subsheaves.

In Section \ref{section:appl} we prove the existence of an effective bound $M_0$ for $M$ which depends only on $X, H, c_1(\widehat\cE)$ and $c_2(\widehat\cE)$ (Theorem \ref{thm:mu-stability2}). In case of elliptic fibrations a similar bound has been proven
for elliptic surfaces (cf.~Theorem 7.4, \cite{FMW99}) and it was expected that such a bound exists
also in higher dimensions. As a result we find that the Fourier-Mukai transform induces an embedding of the relative Jacobian of line bundles of degree $d$ on a curve (in the moduli space of curves $\cM([C],g)$ of arithmetic genus $g$ in the cohomology class $[C]$) into the moduli space of sheaves on $X$ with Chern classes $rn, c_1,ch_2, ch_3$, in particular, this makes the moduli space of $H_{M_0}$-stable spectral sheaves with Chern classes $rn, c_1,ch_2, ch_3$ to a generic torus fibration
over the moduli space $\cM([C],g)$.

In Section \ref{section:case} we study the special case of $r=1$ which was not covered in the previous sections. Here the kernel of the Fourier-Mukai transforms is given by the ideal sheaf of the relative
diagonal immersion $\delta\colon X\hookrightarrow X\times_B X$. We prove that this kernel
gives an equivalence of categories.

The case of trivial fibrations $X=S\times B \to B$ where $S$ is a K3 surface and $B$ is an elliptic curve is covered in Section \ref{sec:trivial}. If we assume additionally that $S$ is reflexive in the sense of \cite{BBH97a} all the topological invariants of the universal bundle are known, which allows for more specific computations.

We are then able to construct many instances of stable sheaves (even with vanishing first Chern class) on K3 fibered Calabi-Yau threefolds out of spectral curves embedded in them.

\subsubsection*{\bf Acknowledgments.} We would like to thank T. Bridgeland for discussions and the anonymous referee for comments and suggestions which helped us to improve the paper. B. Andreas thanks the SFB 647, A1 ``Space-Time-Matter, Analytic and Geometric Structures'' and the Free University Berlin for support where part of this work was carried out. D.~Hern\'andez Ruip\'erez and D.S\'anchez G\'omez are supported by the spanish grant MTM2006-04779 (Ministry for Education and Science) and by grants SA001A07 and GR46 (Junta de Castilla y Le\'on).

\subsubsection*{\bf Conventions} By a variety we understand a scheme locally of finite type over the complex numbers. We denote by $\dbc{X}$ the bounded derived category of complexes of quasi-coherent sheaves with coherent cohomology. An integral functor $\Phi=\fmf{\cplx K}{Y}{X}\colon \dbc{Y} \to \dbc{X}$ is called a Fourier-Mukai  transform when it is an equivalence of categories and its kernel $\cplx K$ reduces to a single sheaf. We denote by $\Phi_t$ the integral functor on the fiber defined by the derived restriction of the
kernel ${\cplx K}$ of $\Phi$. If $x\in X$ is a point, $\mathcal{O}_x$ denotes the skyscraper sheaf of length 1 at $x$, that is, the structure sheaf of
$x$ as a closed subscheme of $X$. By a Calabi-Yau threefold we understand a smooth projective variety of dimension 3 with trivial canonical bundle.

\section{Moduli Spaces and Fourier-Mukai transform\label{section:moduli}}
In this paper we will consider K3 fibrations $X$ which satisfy the following conditions.
We assume that $p\colon X\to B$ is a proper flat morphism of non-singular projective varieties with connected fibers whose generic fiber is a smooth K3 surface and $B$ is an irreducible curve.

Let $H$ be a polarization on $X$. For each $t\in\PP,\,H$ induces a polarization $H_t$ on the fiber $X_t$ of $p$, defined by restricting $H$ to $X_t$.  There exists a relative moduli scheme
$$\hat{p}\colon\mathcal{M}^H(X/\PP)\to\PP$$
whose fiber over a point $t\in\PP$ is the moduli scheme $\mathcal{M}^{H_t}(X_t)$ of stable sheaves on the fiber with respect to $H_t$, and the points of $\mathcal{M}^H(X/\PP)$ correspond to stable sheaves whose scheme-theoretic support is contained in some fiber $X_t$ of $p$ (cf.~\cite{Simp96a}). Here stability means Gieseker stability as considered in \cite[Sect. 1]{Simp96a}. The morphism $\hat p$ takes a sheaf supported on the fiber $p^{-1}(s)$ to the corresponding  point $t\in\PP$.  Let $Y$ be a irreducible component of $\mathcal{M}^H(X/\PP)$ whose points represent sheaves $\cE$,  supported on the fibers $X_t$ which are stable with respect to $H_t$, and have fixed numerical invariants, and hence fixed Hilbert polynomial. Suppose that $Y$ is projective and fine, then there are no properly semistable sheaves and there is a universal sheaf $\cP$ on $Y\times X$ supported on the closed subscheme $Y\times_{\PP} X$ and flat over $Y$. For each point $y\in Y$ the restriction ${\cP}_y$ of $\cP$ to the fiber $\{y\}\times X_t\hookrightarrow X$, with $\hat{p}(y)=t$, is the stable sheaf represented by the point $y$. The universal sheaf $\cP$ defines a relative integral functor $\Phi\colon \dbc{Y}\to \dbc{X}$ between the bounded derived categories of $Y$ and $X$ by letting
$$\Phi(\cplx{E})=\bR{\pi_X}_\ast(\pi_Y^\ast(\cplx{E})\otimes\cP)\,,$$
where $\pi_Y$ and $\pi_X$ denote the projections of $Y\times_{\PP} X$ onto its factors.

When $Y$ and $X$ have the same dimension and $q\colon Y\to\PP$ is equidimensional,  by a result due to Bridgeland and Maciocia \cite[Theorem 1.2 ]{BrM02}, one has that $Y$ is smooth, $q$ is a K3 fibration and the integral functor $\Phi$ is a Fourier-Mukai transform, that is an equivalence of categories.

To construct such components $Y$ of the relative moduli space $\mathcal{M}^H(X/\PP)$, we
assume that there exists a divisor $L$ on $X$ and integer numbers $r>0$, $s$, 
such that there exists a sheaf $\cE$ on a non-singular fiber $i_t\colon X_t\hookrightarrow X$, which is stable with respect to $H_t$ and has Mukai vector $v=(r,L_t,s=\ch_2(\cE)+r)$, where $L_t$ is the restriction of $L$ to $X_t$. 
The Hilbert polynomial of $\cE$ is
\begin{equation}\label{eq:Hilbert polynomial}
P(n)=\frac{1}{2}rH_t^2 n^2+L_t H_t n+(r+s)= rH_t^2\binom{n+1}{2}+ a(L_t) n + (r+s)\,,
\end{equation}
with $a(L_t)=L_t H_t - \frac12 rH_t^2\in \mathbb Z$.
\begin{prop}\label{prop:Y}
If the greatest common divisor of $rH_t^2$, $a(L_t)$, and $s+r$ is equal to 1 and the Mukai vector $v$ satisfies $v^2=0$, then there exists an irreducible component $Y$ of $\cM^H(X/B)$ containing the class of sheaf $\cE$ such that:
\begin{enumerate}
\item $Y$ is a fine projective relative moduli space,
\item $\dim Y=3$,
\item $q\colon Y\to B$ is equidimensional.
\end{enumerate}
\end{prop}
\begin{proof}
Let  $\cM^H(X/B, P(n))$ be the relative moduli space of stable sheaves on the fibers of $X\to B$ with Hilbert polynomial $P(n)$.
The coprimality conditions imply that there are no strictly semistable sheaves on the fibers with Hilbert polynomial $P(n)$, so that $\cM^H(X/B, P(n))\to B$ is projective (see, for instance, \cite[Thm.~1.21]{Simp96a}), and that there exists a relative universal on $\cM^H(X/B, P(n))\times_B X$ (see \cite[Theorem A.6]{Muk87a}, which also holds true in the relative setting).

We take the irreducible component $Y$ of $\cM^H(X/B,P(n))$ containing the class of the sheaf $\cE$ whose existence we have assumed. Thus, there exists a universal sheaf $\cP$ on $Y\times_B X$. 

Let us now derive conditions such that the dimension of $Y$ is equal to 3 and the fibration $q$ is equidimensional.

If $X_{t_0}$ is the fiber supporting $\cE$, then $Y_{t_0}$ is nonempty. Since $c_1(\cE)=L_{t_0}$ is the restriction of a divisor $L$ on $X$, the deformations of $\det \cE$ are unobstructed; by Serre duality the traceless part of the second ext's group is $\Ext^2_0(\cE,\cE)=\Hom_0(\cE,\cE)^\ast=0$, and then the deformations of $\cE$ are unobstructed as well. It follows that $Y_t$ is nonempty for all points $t$ in an open neighborhood $U$ of $t_0$ (cf,~\cite[Thm.~4.5]{Th00a}). The image of $q$ is a closed subset of $\PP$ because $Y$ is projective; moreover it contains $U$ and then $q$ is surjective because $B$ is irreducible. We can assume that $X_t$ is smooth for every point $t\in U$. Then  the fiber $Y_t$ of $q\colon Y\to\PP$ is the moduli space of stable sheaves on the K3 surface $X_t$ with Mukai vector $v$, and, by Mukai \cite{Muk84}, $Y_t$ is a smooth variety of dimension $v^2+2=2$. Thus $\dim q^{-1}(U)=3$ so that $Y$ has dimension 3 as well because it is irreducible.

Finally, $q$ is equidimensional because the irreducibility of $Y$ prevents the dimension of the fibers of $q$ from jumping.
\end{proof}
\begin{rem}\label{rem:smooth locus}
Following \cite[Proposition 4.22]{BBH08}, for every smooth fiber $X_t$, the integral functor $\Phi_t$ is an equivalence, so by \cite[Proposition 3.2]{Huy06}, either $Y_t\simeq X_t$ or  the fiber $Y_t$ is a smooth K3 surface which parametrizes $\mu$-stable vector bundles on $X_t$ with $v=(r,L_t,s)$ and $r>1$.
\end{rem}
As we are interested in constructing sheaves on $X$ we need to know when the Fourier-Mukai
transform $\Phi$ maps sheaves to sheaves. As mentioned in the introduction, for this we employ the notion of WIT-sheaves whose definition we now recall.
\begin{defin}
A sheaf $\cE$ on $Y$ is WIT$_i$, with respect to $\Phi$, if there is a coherent sheaf $\cG$ on $X$ such that $\Phi(\cE)\simeq\cG[-i]$.
\end{defin}

For simplicity we shall use the notation $\Phi^j$ to denote the $j$-th cohomology sheaf of $\Phi$, unless confusion can arise. If $\cE$ is WIT$_i$, then $\Phi^j(\cE)=0$ for all $j\neq i$ and we shall write $\widehat{\cE}$   instead of $\Phi^i(\cE)$. So $\Phi(\cE)\simeq\widehat{\cE}[-i]$. Moreover $\widehat{\cE}$ is WIT$_{2-i}$ with respect to $\widehat\Phi$ and $\widehat{\widehat{\cE\,}}\simeq\cE$.
As usual $\widehat{\Phi}$ denotes the quasi-inverse of $\Phi$ shifted by $[-2]$, that is $\Phi\circ\widehat\Phi\simeq [-2]$ and $\widehat\Phi\circ\Phi\simeq [-2]$.

Finally, let us give a concrete example to illustrate how the above conditions can be solved.
\begin{exa}
One of the K3 fibered Calabi-Yau threefolds $X$ which has been analyzed in great detail using
mirror symmetry \cite{Can97} is obtained by resolving singularities of degree eight hypersurfaces $\widehat{X}\subset {\mathbb P}_{1,1,2,2,2}$. The K\"ahler cone of $X$ is generated by positive linear
combinations of the linear system $H=2l+e$ and $l$ where $e$ is an exceptional divisor coming from
blowing-up a curve of singularities and the linear system $l$ is a pencil of K3 surfaces. The intersection ring has been computed in \cite{Can97} and is given by the relations
$$H\cdot l^2=0\,, \quad l^3=0\,, \quad H^3=8\,, \quad H^2\cdot l=4\,.$$
So let us choose $L=xH+yl$ with $x,y\in {\mathbb Z}$  we then find
$$H_t^2=H^2\cdot l=4\,,\quad L_t\cdot H_t=L\cdot H\cdot l=4x\,, \quad L_t^2=L^2\cdot l=4x^2\,.$$
\end{exa}

\section{Stable Spectral Sheaves\label{section:stable}}

\subsection{Stability and relative semistability\label{ssect:stable}}

In this section we give a criterion of stability of torsion free sheaves with respect to $H+M\mathfrak{f}$ for $M\gg 0$ on a K3 fibration $p\colon X \to B$ as above, assuming that
their quotients are good enough.

\begin{defin} The relative degree of a sheaf $\cF$ on $X$ is the intersection number
$$d(\cF)=c_1(\cF)\cdot H\cdot\f$$
where $\f\in A^1(X)$ denotes the algebraic equivalence class of the fiber of $p\colon X\to B$. The relative slope is defined as
$$\mu_{\mathfrak f}(\cF)=\frac{d(\cF)}{\rk(\cF)}\,,
$$
whenever $\rk(\cF)\neq 0$.
\end{defin}
Note that the restriction of a sheaf $\cF$ on $X$ to a general fiber of $p$ has degree $d(\cF)$. Moreover if $\cF$ is a sheaf on $X$ flat over $\PP$, its relative degree is the degree of the restriction ${\cF}_t$ to any fiber $X_t$ of $p$, so that $\mu_{\mathfrak f}(\cF)=\mu({\cF}_t)$.

\begin{defin} \label{def:inequ} We say that a torsion free sheaf $\cF$ on $X$ has good quotients if for every subsheaf $\cF'$ of $\cF$ with $0<\rk(\cF')<\rk(\cF)$ one has  $\mu_{\mathfrak f}(\cF') < \mu_{\mathfrak f}(\cF)$.
\end{defin}

\begin{prop}\label{prop:mu-stability}
Let $\cF$ be a torsion free sheaf on $X$ with good quotients. There exists a non-negative integer $M_0$ depending on $\cF$ and $H$ such that $\cF$ is $\mu$-stable with respect to $H+M\f$ for all $M\geq M_0$.\end{prop}
\begin{proof}
Let $\cF'$ be a subsheaf of $\cF$ such that $0<\rk(\cF')<\rk(\cF)$. The $\mu$-stability condition which needs to be solved for $M$ is:
$$\big[\mu_H(\cF')-\mu_H(\cF)\big]+2M\big[\mu(\cF'_t)-\mu(\cF_t)\big]<0\,.$$
Since $\cF$ has good quotients, we have  $\mu(\cF'_t)<\mu(\cF_t)$, where subscripts denote restriction to a general fiber of $p$.
As the set of the slopes $\mu_H(\cF')$, where $\cF'$ ranges over all nonzero subsheaves of $\cF$,  is bounded by the existence of the maximal destabilizing subsheaf of $\cF$ and $\mu_H(\cF)$ is fixed, $\mu_H(\cF')-\mu_H(\cF)$ is bounded as well. Thus, for $M\geq M_0\geq0$ with
$$M_0:= \displaystyle\max_{\cF'\subset \cF}\{\mu_H(\cF')-\mu_H(\cF)\}/2\,,$$
the sheaf $\cF$ is $\mu$-stable with respect to $H+M\f$.
\end{proof}

\subsection{Stability of spectral bundles\label{sec:sspectral}}
Let $C\stackrel{i}\hookrightarrow Y$ be a reduced curve of genus $g$ such that $C\to\PP$ is a flat covering of degree $n$, and $\calL$ a line bundle on $C$. The covering $C\to\PP$ could have at most a finite number of ramification points and $C$ intersects the generic fiber in $n$ different points.
The line bundle $\calL$ defines a sheaf $\cE=i_\ast\calL$ on $Y$, supported on $C$.

\begin{prop}\label{prop:WIT0}
$\cE$ is WIT$_0$ with respect to the Fourier-Mukai transform $\Phi$, and the transformed sheaf $\what\cE=\Phi(\cE)$ is flat over the base $B$.
\end{prop}
\begin{proof}
The restriction of $\cE$ to each fiber has 0-dimensional support and since the cover $C\to\PP$ is flat, $\cE$ is flat as well. Then, by \cite[Corollary 1.8]{BBH08}, $\cE$ is WIT$_0$ with respect to $\Phi$, $\widehat{\cE}$ is flat over $B$ and $(\widehat{\cE}\,)_t\simeq\widehat{{\cE}_t}$ for every point $t\in B$.
\end{proof}
If $C_t$ denotes the restriction of $C$ to any fiber $Y_t$, then $C_t=\bigsqcup_{y\in C_t} Z_{y}$ where $Z_{y}\hookrightarrow C$ is a zero-cycle of length $n_{i(y)}$ supported on $y$ and such that $\sum n_{i(y)}=n$. If $C\to B$ is not ramified at $y$ one has $n_{i(y)}=1$, so $n_{i(y)}>1$ at only finitely many points $y\in C$.
Therefore $\cO_{C_t}$ has a filtration with factors isomorphic to $\cO_y$, for every $y\in C_t$. 
By base-change $\widehat{{\cE}_t}=\Phi_t(i_\ast\cO_{C_t})$ and we have also a filtration for $\widehat{{\cE}_t}$ whose factors are isomorphic to $\cP_y$, for every $y\in C_t$. In particular, if $Y_t$ is a generic fiber, then $C$ intersects $Y_t$ at $n$ different points $\{y_1,\ldots,y_n\}$ and $\widehat{{\cE}_t}=\cP_{y_1}\oplus\cdots\oplus\cP_{y_n}$.
\begin{prop}\label{prop:torsion free}
$\what\cE$ is a torsion free sheaf on $X$. Moreover, the restriction $\widehat{{\cE}_t}$ of the transformed sheaf $\what\cE$ to any fiber $X_t$ is semistable of degree $n L\cdot H\cdot \f$ with respect to $H_t$.
\end{prop}
\begin{proof}
The restriction $\what\cE_t$ of $\what\cE$ to any fiber has a filtration with quotients of the form $\cP_y$. Since the support of $\cP_y$ is $X_t$ and $\cP_y$ is semistable of degree $L\cdot H\cdot \f$ with respect to $H_t$, for all $y\in C$, both statements of the proposition hold true.
\end{proof}
An alternative description of spectral sheaves is as follows. Let us consider the diagram
\begin{equation*}
\xymatrix{C\times_B X\ar[d]_{p_C}\ar@{^{(}->}[r]^j & Y\times_B X\ar[r]^(0.6){\pi_X}\ar[d]_{\pi_Y} & X\ar[d]^p\\
C\ar@{^{(}->}[r]^i & Y\ar[r]^q & B}
\end{equation*}
where $p_C$ is a flat morphism and we denote $p_X=\pi_X\circ j$. By base change the spectral sheaf $\what\cE$ can be described as

$$\Phi^{\cP}_{Y\to X}(i_\ast\calL)={p_X}_\ast(p_C^\ast\calL\otimes \bL j^\ast\cP)\,.$$

For each point $y\in C$ with $q(y)=t$, the derived inverse image $\bL i_t^\ast\bL j^\ast \cP$ is isomorphic to the sheaf $\cP_y$ where as usual $i_t\colon \{y\}\times X_t\hra C\times_B X$ denotes the immersion of the fiber $p_C^{-1}(y)$. By \cite[Proposici\'on 1.11]{BBH08} $\bL j^\ast \cP$ is a coherent sheaf $\rest{\cP}{C\times_B X}$ on $C\times_B X$, flat over $C$, and such that $i_t^\ast\rest{\cP}{C\times_B X}\simeq\cP_y$ for every $y\in C$. Then
$$\what\cE={p_X}_\ast(p_C^\ast\calL\otimes\rest{\cP}{C\times_B X})\,.$$

\begin{rem}
Note that for elliptic fibrations one shows that any torsion free
sheaf $\cF$ that is semistable of degree zero on the fiber, has
the $WIT_1$ property with respect to the inverse Fourier-Mukai transform
and so gives rise to a single sheaf. For K3 fibrations there is no obvious
condition which guarantees that the inverse Fourier-Mukai transform of a
sheaf $\cF$ (semistable on the fibers) is a single sheaf. One problem is
that the derived dual of $\cP$ may fail to be a single sheaf, in contrast to the elliptic
fibrations where the derived dual of the Poincar\'e sheaf is again a single sheaf. Therefore
in this paper we will only construct sheaves out of spectral data.
\end{rem}

We now show that there exists an integer number $M_0$, such that $\what{\cE}$ is a $\mu$-stable torsion free sheaf on $X$ with respect to $H+M\f$ for all $M\geq M_0$.

\begin{prop}\label{prop:inequ}
For every torsion free quotient $\cG\neq \what\cE$ of $\widehat\cE$ the inequality $\mu(\cG_t)>\mu(\widehat{\cE_t})$
holds true for the generic fiber if and only if the spectral cover $C$ is reduced and irreducible.
\end{prop}

\begin{proof}
If $C$ is not reduced and irreducible, then there is a reduced irreducible proper closed subcurve $C'\subset C$ which is a flat covering of $B$ of degree $0<n'<n$. Thus the surjection $\cO_C\to\cO_{C'}$ leads to a surjection $\calL\to\rest{\calL}{C'}$. Applying the Fourier-Mukai transform we obtain a torsion free quotient $\widehat{\cE}\to \cG$ of rank $rn'$ and relative degree $n'\cdot L\cdot H\cdot \f$, which implies $\mu(\cG_t)=\mu(\widehat{\cE_t})$.

Conversely, let $\what\cE \to \cG \to 0$ be torsion free quotient. On the generic fiber we have $\what\cE_t \to \cG_t \to 0$. Since $\what\cE_t$ is semistable, $\mu(\cG_t)\geq \mu(\what\cE_t)$. Moreover, $\what\cE_t=\cP_{y_1}\oplus\dots\oplus \cP_{y_n}$ for $n$ different points $y_1,\dots,y_n$ of $Y_t$, and then, if one has $\mu(\cG_t)= \mu(\what\cE_t)$, the sheaf $\cG_t$ is $\mu$-semistable and the quotients of its Jordan-Holder factors are some of the sheaves $\cP_{y_1},\dots\cP_{y_n}$. 

Thus $\cG_t$ is WIT$_2$, and  then $\cG$ is isomorphic to $\Phi^0{\widehat\Phi^2}(\cG)$ on the inverse image to $X$ of an open subset of $B$.
Applying the inverse Fourier-Mukai transform $\widehat{\Phi}$ to the exact sequence $\what\cE \to \cG \to 0$ one has
$$0\to\cK\to i_\ast\calL\to\widehat\Phi^2(\cG)\to 0\,.$$
Since C is irreducible we have $\cK\simeq i_\ast\calL'$ where $\calL'=\calL\otimes{\cI}_Z$ for some zero-dimensional closed subscheme $Z$ of $C$, and $\widehat\Phi^2(\cG)=i_\ast\cO_Z$. Then  $\Phi^0{\widehat\Phi^2}(\cG)$ is concentrated on a finite number of fibers, which is not possible since $\cG$ is torsion free. Thus we get a contradiction and so $\mu(\cG_t)>\mu(\widehat{\cE_t})$.
\end{proof}

\begin{cor} \label{cor:inequ} For every subsheaf $\cF'$ of $\cF=\widehat{\cE}$ with $0\neq\cF'\neq \cF$, one has $\mu_{\mathfrak f}(\cF') \leq \mu_{\mathfrak f}(\cF)$. Moreover,  $\widehat{\cE}$ has good quotients in the sense of Definition \ref{def:inequ} if and only if the spectral curve $C$ if irreducible and reduced.
\end{cor}
\begin{proof} Let us consider the exact sequence
$$
0 \to \cT \to \cF/\cF' \to \bar\cF \to 0
$$
where $\bar \cF$ is torsion free. We have exact sequences
$$
0\to \cF'' \to \cF \to \bar\cF\to 0\,, \qquad 0 \to \cF' \to\cF'' \to \cT \to 0\,,
$$
Since $\cT$ is a torsion sheaf we get $d(\cF')=c_1(\cF')H\f \leq c_1(\cF'')H\f=d(\cF'')$ and then $\mu_{\mathfrak f}(\cF') \leq \mu_{\mathfrak f}(\cF'')$ because $\cF'$ and $\cF''$ have the same rank. If $C$ is irreducible and reduced and $\rk(\cF')=\rk(\cF)$, one has $\bar\cF=0$ so that $\mu_{\mathfrak f}(\cF') \leq \mu_{\mathfrak f}(\cF)$. Otherwise, one has $\mu_{\mathfrak f}(\cF'')< \mu_{\mathfrak f}(\cF)$ by Proposition \ref{prop:inequ}; hence $\mu_{\mathfrak f}(\cF') < \mu_{\mathfrak f}(\cF)$. Proposition \ref{prop:inequ} also implies that if $\cF$ has good quotients,  the spectral curve $C$ if irreducible and reduced.
\end{proof}

Propositions \ref{prop:mu-stability}, \ref{prop:torsion free} and Corollary \ref{cor:inequ} give the sought after stability result:
\begin{thm}\label{thm:mu-stability}
Assume that the spectral cover $C$ is irreducible and reduced. If $\calL$ is a line bundle on $C$ and $\cE=i_*\calL$,
there exists a non-negative integer $M_0$ depending on $\widehat\cE$, such that the spectral sheaf $\widehat\cE$ is $\mu$-stable with respect to $H+M\f$ for all $M\geq M_0\,$.
\qed\end{thm}

\subsection{Vector Bundles in the case $r>1$\label{subsec:vbs}}
Here we assume $r>1$ and we study the question when $\widehat{\cE}$ is a locally free sheaf .

\begin{prop} \label{prop:cover} If for any singular fiber $Y_t$ there exists a point $y\in Y_t$ such that $\cP_y$ is locally free on $X_t$, then one can choose the covering $C$ in such a way that  $\widehat{\cE}$ is locally free.
\end{prop}
\begin{proof} Note first that since $\widehat{\cE}$ is flat over $\PP$, it is locally free if and only if  $\widehat{\cE}_t$ is locally free for all $t\in\PP$.
Assume now that the condition of the statement holds true. Since locally freeness is an open property (\cite[Lemma 2.18]{HL97}), for every singular fiber $Y_t$ there exists an open subset $U_t\subseteq Y_t$ such that $\cP_y$ is locally free on $X_t$ for every $y\in Y_t$. Moreover, when $Y_t$ is a smooth fiber (and then $X_t$ is a smooth fiber as well), we already know (see Remark \ref{rem:smooth locus}) that $\widehat{\cE}_t$ is locally free. Then, the set $U$ of all points $y\in Y$ such that $\cP_y$ is locally free on the fiber of $y$ is an open subset, which is the complement of a finite number of curves. Thus we can choose $C$ contained in $U$, and then $\widehat{\cE}$ is locally free.
\end{proof}

For elliptically fibered Calabi-Yau threefolds (see \cite{FMW99}), $U_t$ is the smooth locus of $Y_t$ because $\cP_y$ is nothing but the ideal sheaf of the point $y$, so it is locally free if and only if $y$ is a non singular point of $Y_t$.

The hypothesis of Proposition \ref{prop:cover} is fulfilled when the following conditions are fulfilled:
\begin{enumerate}
\item
The fibers of $X\to B$ have at worst ordinary double points as singularities,
\item each fiber $Y_t$ parametrizes $\mu$-stable sheaves,  and 
\item the polarization on $X$ is chosen in such a way that a sheaf, with the Chern classes fixed to define $Y$, is stable if and only if it is $\mu$-stable on the generic fiber.
\end{enumerate}
Actually, Thomas shows \cite[theorem 4.15]{Th00a}  that when both these conditions and the assumptions of Section \ref{section:moduli} hold true, then $Y_t$ has at worst ordinary double point singularities. Moreover singular points represent reflexive non locally free sheaves on $X_t$ and all other points correspond to locally free sheaves.

\begin{rem}\label{rem:trivial}
For the trivial fibration $S\times B\to B$, where $S$ is a K3 surface,  all fibers are smooth. In this case the spectral sheaves $\widehat{\cE}$ are vector bundles. Moreover if the K3 surface is reflexive, in the sense of \cite{BBH97a}, then we have a good description of the universal sheaf $\cP$ which allows the computation of the whole Chern character of $\widehat{\cE}$.
\end{rem}

\begin{rem}
Under the assumptions of \cite{Th00a}, if $B={\mathbb P}^1$ and the cover C has degree equal to one, that is C can be thought of as a closed immersion of ${\mathbb P}^1$ into $Y$, then the spectral sheaf $\widehat\cE$ corresponding to $\cO_C$ is the sheaf ${\tilde{\sigma}}^*T$ in the notation of \cite[Theorem 4.18]{Th00a}, which is locally free.
\end{rem}

\section{$\mu$-Stable Extensions\label{section:ext}}

In this section we construct new $\mu$-stable locally free sheaves $\cF$ as non trivial extensions
$$
0 \to \cE \to\cF \to \cG\to 0 $$
of known ones $\cE$ and $\cG$.
Assume that one has
\begin{enumerate}
\item $\cE$ and $\cG$ have good quotients (Definition \ref{def:inequ}), 
\item the restrictions of $\cE$ and $\cG$ to the general fiber are $\mu$-semistable of degree 0 with respect to $H_t$.
\end{enumerate}

We then know that $\cE$ and $\cG$ are $\mu$-stable with respect to $\bar H=H+M_0 \mathfrak f$ for some $M_0\gg 0$ by Proposition \ref{prop:mu-stability}. Examples will be provided when $\cE$ and $\cG$ are spectral of degree zero on fibers, since they have good quotients  by Corollary \ref{cor:inequ},  and are semistable on the generic fiber. 

We make, in addition, two assumptions which are necessary conditions for $\cF$ to be $\mu$-stable
\begin{enumerate}
\item The extension is not trivial, that is, $\Ext^1(\cG,\cE)\neq 0$,
\item $\mu(\cE)<\mu(\cF)$.
\end{enumerate}

\begin{prop}\label{prop:exts} Let $\cF$ be a non-split extension and assume that  $(ii)$ is satisfied and that
\begin{equation}\label{bound}
\mu_H(\cG)< \mu_H(\cE)+\frac{\rk(\cF)}{\rk(\cE)\rk(\cG)}\,.
\end{equation}
Then  $\cF$ is $\mu$-stable with respect to $H+M\f$ for $M\geq M_0$ sufficiently large depending only on $\cF$.
\end{prop}
\begin{proof} By Proposition \ref{prop:mu-stability}, it is enough to prove that $\cF$ has good quotients, that is, one has $d(\cF')<0$ for any subsheaf $\cF'$ of $\cF$ such that $0\neq \cF'\neq \cF$ and $\cF/\cF'$ is torsion free.

Let us consider the diagram of exact sequences
$$
\xymatrix@R=9pt{
& 0 & 0 & 0 & \\
0 \ar[r] & \cE/\cE' \ar[r]\ar[u] &  \cF/\cF' \ar[r]\ar[u] &  \cG/\cG' \ar[r]\ar[u] & 0\\
 0 \ar[r] & \cE \ar[r]^{i}\ar[u] &  \cF \ar[r]^j \ar[u] &  \cG \ar[r]\ar[u] & 0\\
 0 \ar[r] & \cE' \ar[r] \ar[u] &  \cF' \ar[r] \ar[u] &  \cG' \ar[r]\ar[u] & 0\\
  & 0 \ar[u]& 0\ar[u] & 0 \ar[u]&
}
$$
with $\cE'=i^{-1} \cF'$ and $\cG'=j(\cF')$. Since $d(\cF')=d(\cE')+d(\cG')$ and one has $d(\cE')\leq 0$ and $d(\cG')\leq 0$ by the $\mu$-semistability of $\cE_t$ and $\cG_t$ (for a generic fiber), we have to prove that either $d(\cE')<0$ or $d(\cG')<0$

Assume first that $d(\cE')=0$. Since $\cE$ has good quotients, one has that $\rk(\cE')=\rk(\cE)$ unless $\cE'=0$. Thus, if  $\cE'\neq 0$, the sheaf $\cE/\cE'$ is a torsion subsheaf of $\cF/\cF'$, so that it is zero. It follows that $\cG/\cG'\simeq \cF/\cF'$ is a nonzero torsion free sheaf, which forces $d(\cG')<0$ since $\cG$ has good quotients.

We are then left with the case $\cE'=0$. If $d(\cG')<0$ we conclude; we have then to exclude the case $d(\cG')=0$. Notice that since $\cG$ has good quotients, $d(\cG')=0$ implies  that $\rk(\cG')=\rk(\cG)$. Let us write $D=c_1(\cG/\cG')$; since $\cG/\cG'$ is a torsion sheaf, one has $D\cdot H^2\geq 0$ and $D\cdot H^2=0$ if and only if $\cG/\cG'$ is supported in codimension $\geq 2$.

If $D\cdot H^2>0$, solving the slope condition $\mu(\cG')=\mu(\cF')<\mu(\cF)$ implies the inequality
$$\mu(\cG)- \frac{\rk(\cF)}{\rk(\cG) \rk(\cE)} D\cdot H^2 < \mu(\cE),$$
so that  it suffices to impose \eqref{bound}.
Let us prove that the only remaining case, $D\cdot H^2=0$, cannot possibly happen. We first note that $\cF$ does not admit a destabilizing subsheaf $\cF'=\cG'$ with $\rk(\cG')=\rk(\cG)$ if the map $f\colon \Ext^1(\cG,\cE)\to \Ext^1(\cG',\cE)$ is injective (see \cite{ACu07}, Lemma 2.3). To see this in our case, we apply $\Hom(\ , \cE)$ to the exact sequence
$$0\to \cG'\to \cG\to \cT:=\cG/\cG'\to 0\,,$$
to obtain
$$\Ext^1(\cT,\cE)\to \Ext^1(\cG,\cE)\to \Ext^1(\cG',\cE)\,.$$
Since $\cT$ is supported in codimension $\geq 2$, by duality $\Ext^1(\cT,\cE)=\Ext^2(\cE,\cT\otimes \cO_X(K_X))^*=\Ext^2(\cE,\cT)^*=H^2(X, \cE^*\otimes \cT)^*=0$, so that $f$ is injective, and this finishes the proof.
\end{proof}

\section{Application to moduli spaces\label{section:appl}}

In this section we will assume that $p\colon X\to B$ is a K3 fibered Calabi-Yau threefold. As discussed in Section \ref{section:moduli}, $Y$ is also a smooth threefold. Moreover, $Y$ is Calabi-Yau because one has $\dbc{X}\simeq\dbc{Y}$.

\subsection{Topological Invariants\label{topinv}}

Since $Y$ is a smooth variety, for every object in $\dbc{Y}$ the Chern character is well defined. By the singular Grothendieck-Riemann-Roch theorem \cite[Corollary 18.3.1]{Fu98}, the Chern characters of the Fourier-Mukai transform $\Phi(\cplx{E})$ of an object $\cplx{E}$ in $\dbc{Y}$ can be computed in terms of  those of $\cplx{E}$ as follows
$$
\ch(\Phi(\cplx E))={\pi_X}_\ast\Big(\pi_Y^\ast\big(\ch(\cplx E)\cdot \Td(Y/\PP)\big)\cdot\ch(\cP)\Big)
$$
whenever $\ch(\cP)$ exists. Also note that since $p\colon Y\to\PP$ is a local complete intersection morphism, there exists a virtual relative tangent bundle in the K-group $K^\bullet(Y)$. The following is a straightforward computation
\begin{prop}
The Todd class of the (virtual) tangent bundle $T_{Y/\PP}$ is
$$\Td(Y/\PP)=1-\alpha{\hat{\f}\,}+\frac{c_2(Y)}{12}-2\alpha\varpi_Y\,,$$
where $\alpha=\Td_1(B)$, $\hat\f\in A^1(Y)$ denotes the algebraic equivalence class of the fiber of $q\colon Y\to B$, $\varpi_Y$ is the fundamental class of $Y$, and, as usual, we have identified $H^6(Y,\bZ)$ with $\bZ$ by integrating over the fundamental class of $Y$.
\end{prop}
We now take a reduced spectral curve $C\stackrel{i}\hookrightarrow Y$ of genus $g$ such that $C\to\PP$ is a flat covering of degree $n$, and $\calL$ a line bundle of degree $d$ on $C$. We consider the sheaf   $\cE=i_\ast\calL$ on $Y$, supported on $C$. As we have already seen $\cE$ is WIT$_0$, with respect to $\Phi$, and gives then rise to a transform $\widehat{\cE}$ on $X$. Since on a generic fiber $\widehat{\cE_t}\simeq\cP_{y_1}\oplus\cdots\oplus\cP_{y_n}$ we should get
$$\rk(\widehat{\cE})=n r\quad\text{ and }\quad c_1(\widehat\cE)=nL+k\f\quad\text{ with $k$ an integer number}\,.$$

Note that we cannot obtain sheaves with rank equal to the degree of the cover as in the spectral cover construction on elliptic fibrations where the universal sheaf parametrizes sheaves with rank equal to one.

The Chern character of $\cE$ is given by
$$\ch(\cE)=i_\ast(\ch\calL\cdot (1-\frac{1}{2} K_C))\,,$$
where $K_C$ is a canonical bundle of $C$.
Thus the Chern characters of $\cE$ are
\begin{align*}
{\ch}_0(\cE)&=0\,, \quad\quad {\ch}_1(\cE)=0\,,\quad\quad
{\ch}_2(\cE)=i_\ast(1)=[C]\,, \\
{\ch}_3(\cE)&=i_\ast(-\frac{1}{2}K_C+c_1(\calL))=(1-g+d)\varpi_Y\,.
\end{align*}
The Riemann-Roch theorem on $C$ implies
$\chi:=\chi(\calL)=1-g+d$.

When the universal sheaf $\cP$ has  finite homological dimension as a sheaf on $Y\times_{\PP} X$, the Chern character exists and to compute it we should need an explicit expression of $\cP$. Let us denote $\gamma^i=\ch_i(\cP)$, then $\gamma^0=r$.

We find for $\ch_i(\widehat\cE)$ the following expressions
\begin{equation}
\label{e:chern}
\begin{aligned}
\rk(\widehat{\cE})&={\pi_X}_\ast(\pi_Y^\ast [C]\gamma^0)=r[C]\cdot\hat\f=rn\\
{c}_1(\widehat{\cE})&={\pi_X}_\ast(\pi_Y^\ast [C]\gamma^1)+(\chi-\alpha n)r\f\\
\ch_2(\widehat{\cE})&={\pi_X}_\ast(\pi_Y^\ast [C]\gamma^2)+(\chi-\alpha n)L\f\\
\ch_3(\widehat{\cE})&={\pi_X}_\ast(\pi_Y^\ast [C]\gamma^3)+(\chi-\alpha n)(s-r).
\end{aligned}
\end{equation}
Note that by Hurwitz's formula $\chi-\alpha n=d-\frac{R}{2}$ where $R$ denotes the ramification of $C\to B$.

Since $\cP$ is of finite homological dimension, the determinant bundle $\det\cP$ is defined, and the first Chern character of $\cP$ can be computed as follows. We consider the line bundle $\cM=\det\cP\otimes \pi_X^\ast\cO_{X}(-L)$ on $Y\times_{\PP} X$; for every point $y\in Y$  one has ${\pi_Y}^{-1}(y)\simeq \{y\}\times X_t$ with $t=q(y)$, and by restricting to the fiber ${\pi_Y}^{-1}(y)$ we obtain
$$\cM_{|_{X_t}}\simeq\det\cP_{|_{X_t}}\otimes {\pi_X}^\ast\cO_{X}(-L)_{|_{X_t}}\simeq\det{\cP_t}_{|_{X_t}}\otimes \cO_{X_t}(-L_t)\simeq\cO_{X_t}
$$
Thus $\cM\simeq{\pi_Y}^\ast\cQ$ for a line bundle $\cQ$ on $Y$. 
\begin{defin}
We define the divisor class $Q$ such that $\cQ:=\cO_{Y}(-Q)$.\end{defin}
Then one has
$$\det\cP\simeq{\pi_X^\ast}\cO_{X}(L)\otimes{\pi_Y^\ast}\cO_{Y}(-Q)$$
and $c_1(\what\cE\,)$ can be expressed as
$$c_1(\widehat{\cE})=nL+\big((\chi-\alpha n)r-[C]\cdot Q\big)\f$$
The expression of $c_1(\widehat{\cE})$ shows that we cannot twist $\widehat{\cE}$ by a line bundle
to get vanishing first Chern class of $\widehat\cE$. One cannot choose $L=rL'$ as this would be in contradiction to find a simultaneous solution to $\gcd(rH_t^2,a(L_t),r+s)=1$ and $v^2=0$.

If $\cP$ is not of finite homological dimension, we can only consider the push-forward $j_\ast\cP$ where $j\colon Y\times_B X \hookrightarrow Y\times X$, which always has Chern classes because $Y\times X$ is smooth. The relative integral functor with kernel $\cP$ is equal to the absolute integral functor with kernel $j_\ast\cP$ and one can also compute the Chern characters of $\what\cE$ in terms of $\ch(j_\ast\cP)$ and the Todd class of $Y$. 

\subsection{Moduli Spaces}
Although we can use Theorem \ref{thm:mu-stability} to produce $\mu$-stable torsion free sheaves on $X$, that result cannot be used directly in order to produce isomorphisms between moduli spaces, because the polarization $H=H_0+M\f$ depends on the sheaf $\what{\cE}$. Our next aim is to prove that we can find $M$ an effective bound for $M$, which only depends on the topological invariants of $\what{\cE}$.

Let $H_0$ a smooth ample divisor on $X$, which exists by  Bertini's theorem (c.f. \cite[Theorem 8.18]{Hart77}). We put $H_M=H_0+M\f$ for $M\geq 0$. The following result is a straightforward consequence of \cite[5.3.5 and 5.3.6]{HL97}.
\begin{lem}\label{lem:suitable polarization}
Let $D$ be a divisor on $X$ and $k>0$ an integer number such that $0>D^2H_0\geq -k$. Then either $D H_0\f=0$ or $(D H_0\f)(D H_0 H_M)>0$ for every $M\geq \frac{k}{2}(H_0^2\f)$.
\end{lem}

Let $\cF$ be a coherent sheaf on $X$ with Chern classes $c_1$ and $c_2$  and rank $m$. The discriminant of $\cF$ by definition is the characteristic class
$$
B(\cF)=2m c_2-(m-1){ c_1}^2\,.
$$
The Bogomolov inequality states that
if $\cF$ is torsion free and $\mu$-semistable with respect to a some polarization $H$, then $B(\cF)H\geq 0$ (\cite[Theorem 7.3.1.]{HL97}).
\begin{thm}\label{thm:mu-stability2}
Let $H_0$ be a smooth ample divisor on $X$. Then $\what{\cE}$ is $\mu$-stable with respect to $H_M=H_0+M\f$ for every $M\geq M_0=\frac{r^2n^2}{8}B(\what{\cE})H_0(H_0^2\f)$.
\end{thm}
\begin{proof}
Assume $\what{\cE}$ is not $\mu$-stable with respect to $H_{M_0}$. By Theorem \ref{thm:mu-stability},  $\what{\cE}$ is $\mu$-stable with respect to $H_M$ for $M\gg 0$. There exists a  rational number $M_1\geq M_0$ such that $\what{\cE}$ is strictly $\mu$-semistable with respect to $H_{M_1}$ and $\mu$-stable with respect to $H_M$ for all $M>M_1$. So, there is an exact sequence
$$0\to\cF\to\what{\cE}\to\cG\to 0$$
where $\cF$ and $\cG$ are torsion free and $\mu$-semistable sheaves with respect to $H_{M_1}$. We consider
the divisor $D=\rk(\cF)c_1(\cG)-\rk(\cG)c_1(\cF)$, then one has $DH_0H_{2M_1}=DH_{M_1}^2=0$. Notice that $DH_0$ is not numerically trivial, as divisor on $H_0$, because otherwise $\what{\cE}$ would be strictly $\mu$-semistable for all $M$, then by the Hodge index theorem one has $D^2H_0<0$. On the other hand we have
$$B(\what{\cE})H_0=\frac{rn}{\rk(\cF)}B(\cF)H_0+\frac{rn}{\rk(\cG)}B(\cG)H_0-\frac{1}{\rk(\cF)\rk(\cG)}D^2H_0\,.$$
Since $H_{M_1}^2=H_0H_{2M_1}$, the sheaves $\cF$ and $\cG$ are $(H_0,H_{2M_1})$-semistable in the sense of \cite{Miya87}.  By \cite[Corollary 4.7]{Miya87} we have $B(\cF)H_0\geq 0$ and $B(\cG)H_0\geq 0$. Hence
$$0>D^2H_0\geq -\frac{r^2n^2}{4}B(\what{\cE})H_0\,.$$
Since $M_1\geq M_0$ and $DH_0H_{2M_1}=0$ we get $DH_0\f=0$, by Lemma \ref{lem:suitable polarization}; this contradicts Proposition \ref{prop:inequ}.
\end{proof}

Notice that due to Equation \eqref{e:chern}, if we express the discriminant class $B(\what\cE\,)$ in terms of the topological invariants of $\cE=i_*\calL$, we obtain the formula:
$$
B(\what\cE\,)=n^2 L^2 -  2n ([C]\cdot Q) L\cdot \f -2rn (\pi_{X*}(\pi_Y^*[C]\cdot \gamma^2)\,,
$$
where $\gamma^2=\ch_2(\cP)$, which proves that $B(\what\cE\,)$ depends only on the cohomology class of spectral curve $C$. Then, once we have fixed the original polarization $H_0$,  the number $M_0$ depends also only on $[C]$.

Let us now fix the following data:
\begin{enumerate}
\item Three integer numbers $n\ge 1$, $g\ge 0$ and $d$;
\item the cohomology class $[C]$ of a curve $i\colon C \hookrightarrow Y$ of arithmetic genus $g$ which is flat of degree $n$ over $B$.
\end{enumerate}

We denote by  $\cM([C], g)$ the moduli space of integral flat covers of $B$ of arithmetic genus $g$ in the cohomology class $[C]$ and by $\operatorname{Pic}([C],g,d)$ the relative Jacobian of line bundles of degree $d$ on the curves in  $\cM([C], g)$; there is  then a morphism $\operatorname{Pic}([C],g,d) \to \cM([C],g)$ whose fiber over a point $C' \in  \cM([C],g)$ is the Jacobian $\operatorname{Pic}(C',d)$.
Now as $M_0$ depends only on $[C]$, Theorem \ref{thm:mu-stability2} implies that the Fourier-Mukai transform $\Phi$ induces an immersion of algebraic varieties
$$
\Phi\colon \operatorname{Pic}([C],g,d) \hookrightarrow \cM^{H_{M_0}}(X,rn, c_1,\ch_2,\ch_3)
$$
where $\cM^{H_{M_0}}(X,rn, c_1,\ch_2,\ch_3)$ is the moduli space of sheaves on $X$ with Chern classes $rn$, $c_1$, $\ch_2$ and $\ch_3$ given by the right hand sides of Equation \eqref{e:chern}, and which are stable with respect to $H_{M_0}$. The image is the moduli space $\cS^{H_{M_0}}(X,rn, c_1,\ch_2,\ch_3)$ of stable spectral sheaves with those Chern classes, and there is an isomorphism
$$
\Phi\colon  \operatorname{Pic}([C],g,d) \simeq \cS^{H_{M_0}}(X,rn, c_1,\ch_2,\ch_3)\,.
$$
On then has:
\begin{prop} There is a morphism 
$$
\cS^{H_{M_0}}(X,rn, c_1,\ch_2,\ch_3) \to \cM([C],g)
$$ which makes the moduli space of $H_{M_0}$-stable spectral sheaves with Chern classes $(rn, c_1,\ch_2,\ch_3)$, into a generic abelian fibration over the moduli space $\cM([C],g)$ of curves of arithmetic genus $g$ in the cohomology class $[C]$.
\qed
\end{prop}
\section{The Case $r=1$\label{section:case}}
In this section we will treat the rather special case $r=1$ and assume that $X$ is a K3 fibered Calabi-Yau threefold. Note that by adjunction formula every fiber has also trivial canonical bundle. In this case the universal sheaf of Section \ref{intro} is just the ideal sheaf $\cI_\Delta$ of the relative diagonal immersion $\delta\colon X \hookrightarrow X\times_{\PP} X$.
Since for a singular fiber $X_t$ the ideal sheaf of a point $\cI_x$ may be not stable, we cannot use Bridgeland-Maciocia result (\cite[Theorem 1.2 ]{BrM02}) to show that the ideal sheaf $\cI_\Delta$ defines an autoequivalence of the bounded derived category of X. So, we shall prove this directly. Let us denote by $\Phi$ the relative integral functor defined by $\cI_{\Delta}$. In order to prove $\Phi$ maps $\dbc{X}$ to $\dbc{X}$ we need the following notion, which was introduced in \cite{HLS08}.

\begin{defin}
Let $f\colon Z \to T$ be a morphism of
schemes. An object $\cplx E$ in $\dbc{Z}$ is said to be of
\emph{finite homological dimension over $T$}  if $\cplx E\lotimes
\bL f^\ast \cplx G$ is bounded for any $\cplx{G}$ in $\dbc{T}$.
\end{defin}
Following \cite[Proposition 2.7]{HLS08}, $\Phi$ takes bounded complex to bounded complex if and only if $\cI_{\Delta}$ is, as object in $\dbc{X\times_{B}X}$, of finite homological dimension over the first factor. Taking into account the exact sequence
$$0\to \cI_\Delta\to \cO_{X\times_{B}X}\to \cO_{\Delta}\to 0\,,$$
$\cI_{\Delta}$ is of finite homological dimension over the first factor if and only if  $\cO_{\Delta}$ has the same property, and this is indeed the case because the projection formula for $\delta$ implies that  $\cO_{\Delta}=\delta_\ast\cO_X$. We notice that by symmetry the ideal sheaf $\cI_\Delta$ is an object in $\dbc{X\times_{B}X}$ of finite homological dimension over both factors.

\begin{prop}
The relative integral functor
$$\Phi\colon\dbc{X}\to\dbc{X}$$
defined by the ideal sheaf of the relative diagonal is a Fourier-Mukai transform, that is an equivalence of categories.
\end{prop}
\begin{proof}
Since $\cI_\Delta$ is of finite homological dimension over both factors, $\Phi$ is an equivalence if and only if $\Phi_t$ is an equivalence for all $t\in B$, by Proposition 2.15 of \cite{HLS08}. Moreover, the sheaf $\cI_\Delta$ is flat over $B$, and then,  for every closed point $t\in B$ the derived restriction $\bL j_t^\ast\cI_\Delta$ to the fiber $X_t\times X_t$ is the ideal sheaf $\cI_{\Delta_t}$ of the diagonal embedding $X_t\hookrightarrow X_t\times X_t$. Thus $\Phi_t$ is the integral functor defined by $\cI_{\Delta_t}$. A straightforward computation shows that $\Phi_t\simeq T_{\cO_{X_t}}[-1]$, where $T_{\cO_{X_t}}$ denotes the twisted functor along the object $\cO_{X_t}$, that is the integral functor whose kernel is the cone of the morphism  $\cO_{X_t}^\ast\boxtimes\cO_{X_t}\to\cO_{\Delta_t}$. According to \cite{SeTh01}, $\Phi_t$ is an equivalence whenever $\cO_
{X_t}$ is a spherical object of $\dbc{X_t}$.
By hypotheses $X_t$ is a Gorenstein connected surface with trivial dualizing sheaf, for all $t\in B$. Moreover, since the fibration $p\colon X\to B$ is flat and the generic fiber is a (smooth) K3 surface one has $\chi(X_t,\cO_{X_t})=2$ for all $t\in B$. Thus $\cO_{X_t}$ is a spherical object whence the statement follows.
\end{proof}

Let us then consider a spectral cover $C\to B$ of degree $n$ embedded into $X$ by $i\colon C \hookrightarrow X$ and a line bundle $\calL$ of degree $d$ on $C$. By Proposition \ref{prop:torsion free} and Theorem \ref{thm:mu-stability}, the Fourier-Mukai transform $\what\cE=\Phi(\cE)$, with $\cE=i_*\calL$, is a $\mu$-stable torsion free sheaf.

Let $R=2g_C-2-n(2g_B-2)$ be the ramification of $C\to B$. In this case the invariants of the spectral sheaf $\what\cE$ are given by
\begin{align*}
\ch_0(\what\cE) & = n\\
\ch_1(\what\cE) & = (d-R/2) \mathfrak f \\
\ch_2(\what \cE) & =  - [C]\\
\ch_3(\what \cE) & = (d-R/2-n(g_B-1))\varpi_X
\end{align*}
Then, setting $d= R/2$ we get spectral sheaves of degree 0 and $\ch_3(\what\cE)  = -n(g_B-1)\varpi_X$. 

We can also twist the spectral sheaf to get stable torsion free sheaves with vanishing first Chern class. Set $d= n+ R/2$ and write $\tilde \cE= \what\cE\otimes \cO_X(-\mathfrak f)$. We get:
$$
\ch(\tilde\cE)=(n,0,-[C], -n(g_B-3)\varpi_X)\,.
$$
Note that if we impose that $H^1(X,\cO_X)=0$, the base $B$ is necessarily $\mathbb{P}^1$ and one has $\ch_3(\what\cE)=n$ and $\ch_3(\tilde\cE)=3n$.
Thus, we get a great flexibility for having stable torsion free sheaves with $c_1=0$.

This implies that if a spectral curve $C$ of degree $n$ exists, then there are many nonempty moduli spaces of stable sheaves on $X$.

We know that the spectral sheaves $\what\cE=\Phi(\cE)$ are torsion-free. 
Moreover,
by base change, the restriction of $\what\cE$ to a generic fiber is just $\cI_{x_1}\oplus\cdots\oplus \cI_{x_n}$ where $C\cap\f=\{x_1,\ldots,x_n\}$. Notice that for each point $x$ of a generic fiber the ideal sheaf $\cI_x$ is a $\mu$-stable torsion free sheaf on the fiber, but not locally free. However a straightforward computation shows that the dual $\cI_x^\ast$ is locally free. Actually, the  finish this section we shall see that the dual $\what\cE^\ast$ of the spectral sheaf $\what\cE$ is a vector bundle on $X$. Then, the restriction to the generic fibre of the dual $\what\cE^\ast$ of a spectral sheaf is locally free. Actually one has the following result:

\begin{prop}
$\what\cE^\ast$ is a vector bundle on $X$.
\end{prop}
\begin{proof}
Let us consider the diagram
$$\xymatrix{C\times_{\PP} X\ar@{^{(}->}[r]^j\ar[d]_{p_1} & X\times_{\PP} X \ar[r]^{\pi_2}\ar[d]_{\pi_1} & X\ar[d]^p\\
C\ar@{^{(}->}[r]^i & X\ar[r]^p & {\PP}}$$
where squares are cartesian and vertical morphisms are flat. Let us write $\pi=p\circ i$, $p_C=p\circ i$, and $p_2=\pi_2\circ j$.
A straightforward computation shows that $\what\cE^\ast\simeq ({p_2}_\ast p_1^\ast \calL)^\ast\simeq (p\ast p_{C\ast }\calL)^\ast$, and then it is locally free because $p_C$ is flat.
\end{proof}

The topological invariants of $\what\cE^*$ are given by
$$\ch(\what\cE^\ast)=n+(R/2-d)\f\,.$$

\section{Trivial Fibrations\label{sec:trivial}}
In this section we study the case of a trivial K3 fibration $p\colon X\to B$, that is $X\simeq S\times B$, where $S$ is a reflexive K3 surface and $B$ an elliptic curve $B$.

\begin{defin}
A K3 surface $S$ is reflexive if it carries a polarization $h$ and a divisor $l$ such that
\begin{enumerate}
\item $h^2=2$.
\item $l^2=-12$ and $l\cdot h=0$.
\item $e=l+2h$ is non effective on $S$.
\end{enumerate}
\end{defin}
We recall some well known results about reflexive K3 surfaces. For details we refer to \cite{BBH97a}.
\begin{prop}\label{prop:reflexive K3}
Let us fix the following Mukai vector $v=(2,l,-3)$ on the reflexive K3 surface $S$ and, we denote $\cM_v^h(S)$ the moduli space of stable sheaves on $S$ with Mukai vector $v$. The following holds:
\begin{itemize}
\item [1)] $\cM_v^h(S)$ is a non empty fine moduli space parameterizing $\mu$-stable locally free sheaves on $S$ .
\item [2)] The universal sheaf $\cP_S$ on $\cM_v^h(S)\times S$ induces a Fourier-Mukai transform $\dbc{\cM_v^h(S)}\iso\dbc{X}$
\item [3)] $\cM_v^h(S)$ is a reflexive K3 surface and there exists an isomorphism $\psi\colon\cM_v^h(S)\iso S$ as K3 surfaces.
\item [4)] $[\cE]\in\cM_v^h(S)$ if and only if there exists the exact sequence
$$0\longrightarrow\cO_S\longrightarrow\cE(h)\longrightarrow\cI_x(e)\longrightarrow 0\,,$$
with $\Psi([\cE])=x$.
\end{itemize}
\end{prop}

Following Section \ref{section:moduli} there exists a fine relative moduli space $Y\to B$, fixed $H=h\times B$ and $L=l\times B$, and a universal sheaf $\cP$ on $Y\times_{B} X$ such that for every point $t\in B$ the restriction $\cP_t$ is $P_S$. Moreover $Y$ is a trivial fibration as well, that is $Y\simeq \cM_v^h(S)\times B$.

\begin{prop}\label{prop:sequence key}
There exists a non split exact sequence
$$0\longrightarrow{\pi_1}^\ast\cO_X(E)\longrightarrow
\cN\longrightarrow\cI_\Delta\otimes{\pi_2}^\ast\cO_X(E)\longrightarrow 0
\,,
$$
where $E=L+2H=(l+2h)\times B$.
\end{prop}
\begin{proof}
Let us denote $\pi_i$ the natural projection $X\times_{B} X\to X$ and write $\cF=\cO_X(E)$. For every $t\in B$ we have $H^i(S,\cO_S(e))=0$, for all $i\geq 0$, so $\pm e$ is not effective and $\chi(\cO_S(e))=0$. Then $\bR{\pi_1}_\ast\pi_2^\ast\cF=0$ for all $i\geq 0$. By applying $\pi_1$ to the exact sequence
$$0\longrightarrow\cI_\Delta\otimes\pi_2^\ast\cF\longrightarrow\pi_2^\ast\cF\longrightarrow\cO_\Delta\otimes\pi_2^\ast\cF\longrightarrow 0$$
we get
$\cF\simeq\bR{\pi_1}_\ast(\cO_\Delta\otimes\pi_2^\ast\cF)\simeq \bR{\pi_1}_\ast(\cI_\Delta\otimes\pi_2^\ast\cF)[1]$. By Grothendieck-Verdier duality and adjunction we have
\begin{align*}
\Hom(\cF^\vee[-1],\cF^\vee[-1])\simeq &
\Hom(\pi_1^\ast\cF^\vee[-1],\bR\SHom(\cI_\Delta\otimes\pi^\ast_2\cF,\pi_2^\ast\omega_{X/B}))\\
\simeq & \Hom(\pi_1^\ast\cF^\vee[-1]\otimes\cI_\Delta\otimes\pi^\ast_2\cF,\pi_2^\ast\omega_{X/B})\\
\simeq & \Ext^1(\cI_\Delta\otimes\pi^\ast_2\cF,\pi_1^\ast\cF)
\end{align*}
Then $Id\in \Hom(\cF^\vee[-1],\cF^\vee[-1])$ induces a non split exact sequence
$$0\longrightarrow{\pi_1}^\ast\cO_X(E)\longrightarrow
\cN\longrightarrow\cI_\Delta\otimes{\pi_2}^\ast\cO_X(E)\longrightarrow 0\,.
$$
\end{proof}
The restriction to the fiber of $\pi_1$ gives
$$0\longrightarrow\cO_{S}(e)\longrightarrow
\rest{\cN}{\{x\}\times S}\longrightarrow\cI_x\otimes\cO_{S}(e)\longrightarrow 0$$
By Proposition \ref{prop:reflexive K3} $\rest{\cN}{\{x\}\times S}\simeq \cE(h)$ where $\cE$ is a $\mu$-stable vector bundle, with respect to $h$, on $S$ with Mukai vector $(2,l,-3)$.

Then $\cS:=\cN\otimes\pi_2^\ast\cO_X(-H)$ is a sheaf on $X\times_{B} X$, flat over both factors, and for every point $x\in X$ the sheaf $\rest{\cS}{\{x\}\times S}$ corresponds to a point of $\cM_v^h(S)$. In this way we have a morphism
\begin{align*}
X&\stackrel{\Psi}\longrightarrow Y=\cM_v^h(S)\times B\\
x&\mapsto \rest{\cS}{\{x\}\times S}
\end{align*}
which is an isomorphism as on each fiber is an isomorphism. Then
$$(\Psi\times 1)^\ast\cP\otimes\pi_1^\ast\cU\simeq\cS$$
for some line bundle $\cU$ on $X$.
\begin{thm}\label{thm:reflexive FM}
The integral functor $\Phi^{\cS}\colon\dbc{X}\to\dbc{X}$ defined by $\cS$ is an autoequivalence of $\dbc{X}$, that is, a Fourier-Mukai transform.
\end{thm}
\begin{proof}
By definition $\cS$ is an object in $\dbc{X\times_{B}X}$ which has finite homological dimension over both factors and, for every $t\in B$, the restriction $\cS_t$ defines an integral functor  $\Phi_t\colon\dbc{S}\to\dbc{S}$. Since $\rest{\cS_t}{\{x\}\times S}$ is a $\mu$-stable locally free sheaf with Mukai vector $(2,l,-3)$, the integral functor $\Phi_t$ is an equivalence by Proposition \ref{prop:reflexive K3}. Then the integral functor $\Phi^\cS$ is a Fourier-Mukai transform as a consequence of \cite{HLS08}.
\end{proof}
Taking in account the isomorphism $Y\simeq X$ we can apply the spectral construction of Section \ref{section:stable} directly to $X$. We fix a  cover $i\colon C\hookrightarrow X$ flat over $B$ and $\calL$ a line bundle on $C$. Then $\cE=i_\ast\calL$ is WIT$_0$, with respect to $\Phi^S$, and the Fourier-Mukai transform $\widehat\cE$ is a $\mu$-stable vector bundle on $X$.  We finish this section computing Chern character of $\widehat\cE$. Since $S=\cN\otimes\pi_2^\ast\cO_X(-H)$ one has $\Phi^{\cS}(\cE)\simeq \Phi^{\cN}(\cE)\otimes \cO_X(-H)$. Thus $\cE$ is also WIT$_0$ with respect to $\Phi^{\cN}$ and the Chern character of $\Phi^{\cN}(\cE)$ can be computed using Proposition \ref{prop:sequence key}.
\begin{lem}\label{lem:aux}
Let $\cF$ be the Fourier-Mukai transform of $\cE$ with respect to $\Phi^{\cN}$. The Chern character of $\cF$ is:
$$\Big(2n,nE+(2\chi(\calL)+CE)\f,-2nF+\chi(\calL) E\cdot \mathfrak{f}-C,-3\chi(\calL)-CE\Big)\,,$$
where $F$ denotes the generic fiber of $X\to S$.
\end{lem}
\begin{proof}
By Proposition \ref{prop:sequence key} we have the following exact triangle
$$\bR{\pi_2}_\ast\pi_1^\ast(\cE\otimes\cO_X(E))\to\Phi^{\cN}(\cE)\to\Phi^{\cI_\Delta}(\cE)\otimes\cO_X(E)\,.$$
Then the Chern character of $\cF$ is:
$$\ch(\cF)=\ch\Big(\bR{\pi_2}_\ast\pi_1^\ast(\cE\otimes\cO_X(E))\Big)+
\big(\ch(\bR{\pi_2}_\ast\pi_1^\ast\cE)-\ch(\cE)\big)\cdot\ch \cO_X(E)\,,
$$
that is
\begin{align*}
\ch(\cF)= &
{\pi_2}_\ast\pi_1^\ast\Big(\ch\cE\cdot\ch\cO_X(E)\cdot \Td (X/B)\Big)+\\
+ &{\pi_2}_\ast\pi_1^\ast\Big(\ch(\cE)\cdot \Td (X/B)\Big)\cdot \ch\cO_X(E)+\\
-& \ch\cE\cdot\ch\cO_X(E)\,.
\end{align*}

Since $X=S\times B$ one has $\Td(X/B)=(1,0,2F,0)$ with $F\simeq B$ the generic fiber of $X\to S$.
Moreover
$\ch\cE=(0,0,[C],\chi(\calL))$
and $\ch\cO_X(E)=(1,E,-2F,0)$, and one gets the desired formula.
\end{proof}
As a consequence of Lemma \ref{lem:aux} we have:
\begin{prop}\label{prop:chern}
The Chern characters of the Fourier-Mukai transform $\what\cE$ of $\cE$ with respect to $\Phi=\Phi^{\cS}$ are:
\begin{align*}
\ch_0(\what\cE) & = 2n\\
\ch_1(\what\cE) & = nL+ k \mathfrak f\,, \qquad \text{with} \quad k=2\chi(\calL)+ [C]\cdot E \\
\ch_2(\what\cE) &=\chi(\calL) l-([C]\cdot E)h-4nF-[C]\\
\ch_3(\what \cE) & = -5 \chi(\calL) + H\cdot [C]\,.
\end{align*}
\end{prop}

By Remark \ref{rem:trivial}, the spectral sheaves $\what\cE$ are vector bundles. Moreover, their restrictions to the generic fiber of $S\times B \to B$ are semistable of degree zero by Proposition \ref{prop:torsion free}. If in addition $C$ is irreducible and reduced, $\what\cE$ is stable with respect to $H+M\f$ for all $M\geq M_0$ for an integer $M_0$ depending only on $H$ and $[C]$  (cf.~Theorems \ref{thm:mu-stability} and \ref{thm:mu-stability2}).


\def\cprime{$'$}

\end{document}